\documentclass[reqno, 12pt]{amsart}
\usepackage{a4wide}
\usepackage{amscd}
\usepackage{enumerate}
\input xy
\xyoption{all}

\theoremstyle{plain}

\newtheorem*{thm A}{Theorem~A}
\newtheorem*{thm B}{Theorem~B}
\newtheorem*{thm C}{Theorem~B}
\newtheorem*{pro A}{Proposition~A}
\newtheorem*{pro B}{Proposition~B}
\newtheorem*{lem A}{Lemma~A}
\newtheorem*{lem B}{Lemma~B}
\newtheorem*{lem C}{Lemma~C}
\newtheorem*{lem D}{Lemma~D}

\newtheorem*{rem}{Remark}

\newtheorem*{MT1}{Theorem 1}
\newtheorem*{MT2}{Theorem 2}

\newtheorem{theorem}{Theorem}[section]

\newtheorem{lemma}[theorem]{Lemma}
\newtheorem{proposition}[theorem]{Proposition}
\newtheorem{remark}[theorem]{Remark}

\def \N{\nabla}

\def \Q{Q^m}

\def \x{\xi}
\def \e{\eta}
\def \al{\alpha}
\begin{document}

\title[Reeb parallel structure Jacobi operator]{Real hypersurfaces in the complex quadric with Reeb parallel structure Jacobi operator}

\vspace{0.2in}
\author[H. Lee and Y.J. Suh]{Hyunjin Lee and Young Jin Suh}

\address{\newline
Hyunjin Lee
\newline The Research Institute of Real and Complex Manifolds (RIRCM),
\newline Kyungpook National University,
\newline Daegu 41566, REPUBLIC OF KOREA}
\email{lhjibis@hanmail.net}

\address{\newline
Young Jin Suh
\newline Department of Mathematics \& RIRCM,
\newline Kyungpook National University,
\newline Daegu 41566, REPUBLIC OF KOREA}
\email{yjsuh@knu.ac.kr}


\thanks{* This work was supported by grant Proj. Nos. NRF-2018-R1D1A1B-05040381 and  NRF-2019-R1I1A1A 01050300 from National Research Foundation of Korea.}

\begin{abstract}
In this paper, we first introduce the full express of the Riemannian curvature tensor of a real hypersurface $M$ in complex quadric~$Q^{m}$ from the equation of Gauss. Next we derive a formula for the structure Jacobi operator~$R_{\xi}$ and its derivative under the Levi-Civita connection of $M$. We give a complete classification of Hopf real hypersurfaces with Reeb parallel structure Jacobi operator, $\nabla_{\xi}R_{\xi} =0$, in the complex quadric $Q^{m}$, $m \geq 3$.
\end{abstract}

\maketitle

\section{Introduction}\label{section 1}
\setcounter{equation}{0}
\renewcommand{\theequation}{1.\arabic{equation}}
\vspace{0.13in}
For Hermitian symmetric space of compact type different from the above ones, we can give the example of complex quadric ${Q^m}= SO_{m+2}/SO_mSO_2$, which is a complex hypersurface in the complex projective space ${\mathbb C}P^{m+1}$ (see Romero \cite{R1}, \cite{R2}, Smyth \cite{BS}, Suh \cite{S3}, \cite{S4}). The complex quadric can also be regarded as a kind of real Grassmann manifolds of compact type with rank 2 (see Besse \cite{Bes}, Helgason \cite{He}, and Knap \cite{Kna} ).  Accordingly, the complex quadric ${Q^m}$ admits two important geometric structures, a complex conjugation structure $A$ and a K\"ahler structure $J$, which anti-commute with each other, that is, $AJ=-JA$. Then for $m{\ge}2$ the triple $({Q^m},J,g)$ is a Hermitian symmetric space of compact type with rank $2$ and its maximal sectional curvature is equal to $4$ (see Kobayashi and Nomizu \cite{KO}, Reckziegel \cite{R}).

\par
\vskip 6pt

In addition to the complex structure $J$ there is another distinguished geometric structure on $\Q$, namely a parallel rank two vector bundle ${\mathfrak A}$ which contains an $S^1$-bundle of real structures, that is, complex conjugations $A$ on the tangent spaces of $\Q$. The set is denoted by ${\mathfrak A}_{[z]}=\{A_{{\lambda}\bar z}{\vert}\, {\lambda}\in S^1{\subset}{\mathbb C}\}$, $[z] \in {\Q}$, and it is the set of all complex conjugations defined on $\Q$. Then ${\mathfrak A}_{[z]}$ becomes a parallel rank $2$-subbundle of $\text{End}\ T_{[z]}{\Q}$, $[z] \in {\Q}$. This geometric structure determines a maximal ${\mathfrak A}$-invariant subbundle ${\mathcal Q}$ of the tangent bundle $TM$ of a real hypersurface $M$ in $\Q$.  Here the notion of parallel vector bundle ${\mathfrak A}$ means that $({\bar\nabla}_XA)Y=q(X)JAY$ for any vector fields $X$ and $Y$ on $\Q$,
where  $\bar\nabla$ and ~$q$ denote a connection and a certain $1$-form defined on $T_{[z]}{\Q}$, $[z] \in {\Q}$ respectively (see Smyth~\cite{BS}).

\vskip 6pt

Recall that a nonzero tangent vector $W \in T_{[z]}{\Q}$ is called singular if it is tangent to more than one maximal flat in $\Q$. There are two types of singular tangent vectors for the complex hyperbolic quadric~$\Q$:
\begin{itemize}
\item If there exists a conjugation $A \in {\mathfrak A}$ such that $W \in V(A)=\{X \in T_{[z]}{Q^m}{\vert}\,AX=X\}$, then $W$ is singular. Such a singular tangent vector is called {\it ${\mathfrak A}$-principal}.
\item If there exist a conjugation $A \in {\mathfrak A}$ and orthonormal vectors $Z_{1}$, $Z_{2} \in V(A)$ such that $W/||W|| = (Z_{1}+JZ_{2})/\sqrt{2}$, then $W$ is singular. Such a singular tangent vector is called \emph{${\mathfrak A}$-isotropic},
where $V(A)=\{X \in T_{[z]}{Q^m}{\vert}\, AX=X\}$ and $JV(A)=\{X \in T_{[z]}{Q^m}{\vert} \, AX=-X\}$ are the $(+1)$-eigenspace and $(-1)$-eigenspace for the involution $A$
on $T_{[z]}{Q^m}$, $[z] \in {Q^m}$.
\end{itemize}

\vskip 6pt

On the other hand, Okumura \cite{O} proved that the Reeb flow on a real hypersurface in ${\mathbb C}P^m = SU_{m+1}/S(U_1U_m)$ is isometric if and only if $M$ is an open part of a tube around a totally geodesic ${\mathbb C}P^k$ in ${\mathbb C}P^m$ for some $k \in \{0,\ldots,m-1\}$. For the complex $2$-plane Grassmannian $G_2({\mathbb C}^{m+2})= SU_{m+2}/S(U_2U_m)$ a classification was obtained by Berndt and Suh \cite{BS1}. The Reeb flow on a real hypersurface in $G_2({\mathbb C}^{m+2})$ is isometric if and only if $M$ is an open part of a tube around a totally geodesic $ G_2({\mathbb C}^{m+1})$ in $G_2({\mathbb C}^{m+2})$.   For the complex quadric $Q^m = SO_{m+2}/SO_2SO_m$, Berndt and Suh \cite{BS2} have obtained the following result:
\begin{thm A}\label{Theorem A}
Let $M$ be a real hypersurface in the complex quadric $Q^m$, $m\geq 3$. Then the Reeb flow on $M$ is isometric if and only if $m$ is even, say $m = 2k$, and $M$ is an open part of a tube around a totally geodesic ${\mathbb C}P^k$ in $Q^{2k}$.
\end{thm A}
\par
\vskip 6pt

For the complex hyperbolic space ${\mathbb C}H^m$ a classification was obtained by Montiel and Romero~\cite{MR}. They proved that the Reeb flow on a real hypersurface in ${\mathbb C}H^m$ is isometric if and only if $M$ is an open part of a tube around a totally geodesic ${\mathbb C}H^k$ in ${\mathbb C}H^m$ for some $k \in \{0,\ldots,m-1\}$. The classification problems related to the Reeb parallel shape operator, parallel Ricci tensor,
and harmonic curvature for real hypersurfaces in the complex quadric $Q^m$ were recently given in Suh~\cite{S1}, \cite{S3} and \cite{S4} respectively.

\vskip 6pt

The notion of isometric Reeb flow was introduced by Hutching and Taubes~\cite{HT} and the geometric construction of horospheres in a non-compact manifold of negative curvature was mainly discussed in the book due to Eberlein~\cite{E}.

\vskip 6pt

On the other hand, Jacobi fields along geodesics of a given Riemannian manifold $(\bar M,g)$ satisfy a well known differential equation. This equation naturally inspires the so-called Jacobi operator. That is, if $R$ denotes the curvature operator of $\bar M$, and $X$ is a tangent vector field to $\bar M$, then the Jacobi operator $R_X \in \text{End}(T_{z} {\bar M})$ with respect to $X$ at $z \in \bar M$, defined by $(R_XY)(z)=(R(Y,X)X)(z)$ for any $Y \in T_z {\bar M}$, becomes a self adjoint endomorphism of the tangent bundle $T {\bar M}$ of $\bar M$. Thus, each tangent vector field $X$ to $\bar M$ provides a Jacobi operator $R_X$ with respect to $X$. In particular, for the Reeb vector field~$\xi$, the Jacobi operator $R_{\xi}$ is said to be the {\it structure Jacobi operator}.

\vskip 6pt

Actually, many geometers have considered the fact that a real hypersurface $M$ in K\"{a}hler manifolds has {\it parallel structure Jacobi operator} (or {\it Reeb parallel structure Jacobi operator}, respectively), that is,  $\nabla_{X}R_{\xi}=0$ (or $\nabla_{\xi}R_{\xi}=0$, respectively) for any tangent vector field~$X$ on $M$. Recently Ki, P\'erez, Santos and Suh~\cite{KPSS} have investigated the Reeb parallel structure Jacobi operator in the complex space form $M^m(c)$, $c \neq 0$, and have used it to study some principal curvatures for a tube over a totally geodesic submanifold. In particular, P\'erez, Jeong and Suh~\cite{PJS} have investigated real hypersurfaces $M$ in $G_2({\mathbb C}^{m+2})$ with parallel structure Jacobi operator, that is, ${\nabla}_X{R}_{\xi}=0$ for any tangent vector field $X$ on $M$. Jeong, Suh and Woo \cite{JSW}  and  P\'erez and Santos \cite{PS} have generalized such a notion to the recurrent structure Jacobi operator, that is, $({\nabla}_X{R}_{\xi})Y={\beta}(X){R}_{\xi}Y$ for a certain $1$-form $\beta$ and any vector fields $X,Y$ on $M$ in $G_2({\mathbb C}^{m+2})$ or ${\mathbb C}P^m$. In~\cite{JLS}, Jeong, Lee, and Suh have considered a Hopf real hypersurface with Codazzi type of structure Jacobi operator, $(\nabla_{X}R_{\xi})Y =  (\nabla_{Y}R_{\xi})X$, in $G_{2}(\mathbb C^{m+2})$. Moreover, P\'erez, Santos and Suh~\cite{PSS} have further investigated the property of the Lie $\xi$-parallel structure Jacobi operator in complex projective space ${\mathbb C}P^m$, that is, ${\mathcal L}_{\xi}R_{\xi}=0$.

\vskip 6pt

Motivated by these results, in this paper we want to give a classification of Hopf real hypersurfaces in $Q^{m}$ with non-vanishing geodesic Reeb flow and Reeb parallel structure Jacobi operator, that is, ${\nabla}_{\xi}{R}_{\xi}=0$. Here a real hypersurface $M$ is said to be {\it Hopf} if the Reeb vector field~$\xi$ of~$M$ is principal by the shape operator~$S$, that is, $S\xi=g(S\xi, \xi) \xi = \alpha \xi$. In particular, if the Reeb curvature function $\alpha =g(S\xi, \xi)$ identically vanishes, we say that $M$ has {\it a vanishing geodesic Reeb flow}. Otherwise, $M$ has a {\it non-vanishing geodesic Reeb flow}.

\vskip 6pt

Under these background and motivation, first we prove the following:
\begin{MT1}\label{Main Theorem 1}
There does not exist any Hopf real hypersurface in the complex quadric~$Q^{m}$, $m \geq 3$, with Reeb parallel structure Jacobi operator and $\mathfrak A$-principal singular normal vector field, provided with non-vanishing geodesic Reeb flow.
\end{MT1}

\vskip 6pt

Now let us consider a Hopf real hypersurface with $\mathfrak A$-isotropic singular normal vector field~$N$ in $Q^{m}$. Then by virtue of Theorem~$\rm A$ we can give a complete classification of Hopf real hypersurfaces in $\Q$ with Reeb parallel structure Jacobi operator as follows:
\begin{MT2}\label{Main Theorem 2}
Let $M$ be a Hopf real hypersurface in the complex quadric $\Q$, $m \geq 3$, with Reeb parallel structure Jacobi operator and non-vanishing geodesic Reeb flow. If $M$ has the $\mathfrak A$-isotropic singular normal vector field in $\Q$, then $M$ is locally congruent to a tube around the totally geodesic ${\mathbb C}P^k$ in ${Q^{2k}}$, where $m=2k$, and $r \in (0, \frac{\pi}{4}) \cup (\frac{\pi}{4}, \frac{\pi}{2})$.
\end{MT2}

\vskip 17pt

\section{The complex quadric}\label{section 2}
 \setcounter{equation}{0}
\renewcommand{\theequation}{2.\arabic{equation}}
\vspace{0.13in}

For more background to this section we refer to \cite{K}, \cite{KO}, \cite{R}, \cite{S1}, \cite{S2}, \cite{S4} and \cite{SHw}. The complex quadric $Q^m$ is the complex hypersurface in ${\mathbb C}P^{m+1}$ which is defined by the equation $z_1^2 + \cdots + z_{m+2}^2 = 0$, where $z_1,\cdots,z_{m+2}$ are homogeneous coordinates on ${\mathbb C}P^{m+1}$. We equip $Q^m$ with the Riemannian metric which is induced from the Fubini Study metric on ${\mathbb C}P^{m+1}$ with constant holomorphic sectional curvature~$4$. The K\"{a}hler structure on ${\mathbb C}P^{m+1}$ induces canonically a K\"{a}hler structure $(J,g)$ on the complex quadric. For a nonzero vector $z \in \mathbb C^{m+2}$ we denote by $[z]$ the complex span of $z$, that is, $[z]=\mathbb C z = \{\lambda z\,|\, \lambda \in S^{1} \subset \mathbb C\}$. Note that by definition~$[z]$ is a point in $\mathbb C P^{m+1}$. For each $[z] \in Q^m \subset \mathbb C P^{m+1}$ we identify $T_{[z]}{\mathbb C}P^{m+1}$ with the orthogonal complement ${\mathbb C}^{m+2} \ominus {\mathbb C}z$ of ${\mathbb C}z$ in ${\mathbb C}^{m+2}$ (see Kobayashi and Nomizu \cite{KO}). The tangent space $T_{[z]}Q^m$ can then be identified canonically with the orthogonal complement ${\mathbb C}^{m+2} \ominus ({\mathbb C}z \oplus {\mathbb C}\rho)$ of ${\mathbb C}z \oplus {\mathbb C}\rho$ in ${\mathbb C}^{m+2}$, where $\rho \in \nu_{[z]}Q^m$ is a normal vector of $Q^m$ in ${\mathbb C}P^{m+1}$ at the point $[z]$.

\vskip 6pt

The complex projective space ${\mathbb C}P^{m+1}$ is a Hermitian symmetric space of the special unitary group $SU_{m+2}$, namely ${\mathbb C}P^{m+1} = SU_{m+2}/S(U_{m+1}U_1)$. We denote by $o = [0,\ldots,0,1] \in {\mathbb C}P^{m+1}$ the fixed point of the action of the stabilizer $S(U_{m+1}U_1)$. The special orthogonal group $SO_{m+2} \subset SU_{m+2}$ acts on ${\mathbb C}P^{m+1}$ with cohomogeneity one. The orbit containing $o$ is a totally geodesic real projective space ${\mathbb R}P^{m+1} \subset {\mathbb C}P^{m+1}$. The second singular orbit of this action is the complex quadric $Q^m = SO_{m+2}/SO_mSO_2$. This homogeneous space model leads to the geometric interpretation of the complex quadric $Q^m$ as the Grassmann manifold $G_2^+({\mathbb R}^{m+2})$ of oriented $2$-planes in ${\mathbb R}^{m+2}$. It also gives a model of $Q^m$ as a Hermitian symmetric space of rank $2$. The complex quadric $Q^1$ is isometric to a sphere $S^2$ with constant curvature, and $Q^2$ is isometric to the Riemannian product of two $2$-spheres with constant curvature. For this reason we will assume $m \geq 3$ from now on.

\vskip 6pt

For a unit normal vector $\rho$ of $Q^m$ at a point $[z] \in Q^m$ we denote by $A = A_\rho$ the shape operator of $Q^m$ in ${\mathbb C}P^{m+1}$ with respect to $\rho$. The shape operator is an involution on the tangent space $T_{[z]}Q^m$ and
$$T_{[z]}Q^m = V(A_\rho) \oplus JV(A_\rho),$$
where $V(A_\rho)$ is the $(+1)$-eigenspace and $JV(A_\rho)$ is the $(-1)$-eigenspace of $A_\rho$.  Geometrically this means that the shape operator $A_\rho$ defines a real structure on the complex vector space $T_{[z]}Q^m$, or equivalently, is a complex conjugation on $T_{[z]}Q^m$. Since the real codimension of $Q^m$ in ${\mathbb C}P^{m+1}$ is $2$, this induces an $S^1$-subbundle ${\mathfrak A}$ of the endomorphism bundle ${\rm End}(TQ^m)$ consisting of complex conjugations. There is a geometric interpretation of these conjugations. The complex quadric~$Q^m$ can be viewed as the complexification of the $m$-dimensional sphere~$S^m$. Through each point $[z] \in Q^m$ there exists a one-parameter family of real forms of $Q^m$ which are isometric to the sphere $S^m$. These real forms are congruent to each other under action of the center $SO_2$ of the isotropy subgroup of $SO_{m+2}$ at $[z]$. The isometric reflection of $Q^m$ in such a real form $S^m$ is an isometry, and the differential at $[z]$ of such a reflection is a conjugation on $T_{[z]}Q^m$. In this way the family ${\mathfrak A}$ of conjugations on $T_{[z]}Q^m$ corresponds to the family of real forms $S^m$ of $Q^m$ containing~$[z]$, and the subspaces $V(A)$ in $T_{[z]}Q^m$ correspond to the tangent spaces $T_{[z]}S^m$ of the real forms~$S^m$ of $Q^m$.

\vskip 6pt

The Gauss equation for $Q^m \subset {\mathbb C}P^{m+1}$ implies that the Riemannian curvature tensor $\bar R$ of $Q^m$ can be described in terms of the complex structure $J$ and the complex conjugations $A \in {\mathfrak A}$:
\begin{equation}\label{Riemannian curvature tensor}
\begin{split}
{\bar R}(X,Y)Z &=  g(Y,Z)X - g(X,Z)Y + g(JY,Z)JX - g(JX,Z)JY \\
& \quad  - 2g(JX,Y)JZ  + g(AY,Z)AX \\
& \quad - g(AX,Z)AY + g(JAY,Z)JAX - g(JAX,Z)JAY.
\end{split}
\end{equation}

\vskip 6pt

It is well known that for every unit tangent vector $U \in T_{[z]}Q^m$ there exist a conjugation $A \in {\mathfrak A}$ and orthonormal vectors $Z_{1}$, $Z_{2}\in V(A)$ such that
\begin{equation*}
U = \cos (t) Z_{1} + \sin (t) JZ_{2}
\end{equation*}
for some $t \in [0,\pi/4]$  (see \cite{R}). The singular tangent vectors correspond to the values $t = 0$ and $t = \pi/4$. If $0 < t < \pi/4$ then the unique maximal flat containing $U$ is ${\mathbb R}Z_{1} \oplus {\mathbb R}JZ_{2}$.
Later we will need the eigenvalues and eigenspaces of the Jacobi operator $\bar{R}_U = \bar{R}(\,\cdot\, ,U)U$ for a singular unit tangent vector $U$.
\begin{enumerate}[\rm (1)]
\item {If $U$ is an ${\mathfrak A}$-principal singular unit tangent vector with respect to $A \in {\mathfrak A}$, then the eigenvalues of $\bar{R}_U$ are $0$ and $2$ and the corresponding eigenspaces are ${\mathbb R}U \oplus J(V(A) \ominus {\mathbb R}U)$ and $(V(A) \ominus {\mathbb R}U) \oplus {\mathbb R}JU$, respectively.}
\item {If $U$ is an ${\mathfrak A}$-isotropic singular unit tangent vector with respect to $A \in {\mathfrak A}$ and $X$, $Y \in V(A)$, then the eigenvalues of $\bar{R}_U$ are $0$, $1$ and $4$ and the corresponding eigenspaces are ${\mathbb R}U \oplus {\mathbb C}(JZ_{1}+Z_{2})$, $T_{[z]}Q^m \ominus ({\mathbb C}Z_{1} \oplus {\mathbb C}Z_{2})$ and ${\mathbb R}JU$, respectively.}
\end{enumerate}

\vskip 17pt

\section{Real hypersurfaces in $Q^{m}$}\label{section 3}
\setcounter{equation}{0}
\renewcommand{\theequation}{3.\arabic{equation}}
\vspace{0.13in}

Let $M$ be a  real hypersurface in $Q^m$ and denote by $(\phi,\xi,\eta,g)$ the induced almost contact metric structure. By using the Gauss and Wingarten formulas the left-hand side of \eqref{Riemannian curvature tensor} becomes
\begin{equation*}
\begin{split}
{\bar R}(X,Y)Z  & = R(X,Y)Z -g(SY, Z)SX + g(SX, Z)SY \\
& \quad + \big\{g((\nabla_{X}S)Y, Z)- g((\nabla_{Y}S)X, Z) \big \} N,
\end{split}
\end{equation*}
where $R$ and $S$ denote the Riemannian curvature tensor and the shape operator of $M$ in $Q^m$, respectively. Taking tangent and normal components of \eqref{Riemannian curvature tensor} respectively, we obtain
\begin{equation}\label{eq: 2.1}
\begin{split}
& g(R(X,Y)Z, W) - g(SY,Z)g(SX,W) + g(SX,Z)g(SY,W) \\
& =  g(Y,Z)g(X,W) - g(X,Z)g(Y,W) + g(JY,Z)g(JX,W) \\
& \ \ - g(JX,Z)g(JY,W) - 2g(JX,Y)g(JZ, W) + g(AY,Z)g(AX, W) \\
& \ \  - g(AX,Z)g(AY, W)+ g(JAY,Z)g(JAX,W)- g(JAX,Z)g(JAY,W), \\
\end{split}
\end{equation}
and
\begin{equation}\label{codazzi equation}
\begin{split}
& g((\nabla_{X}S)Y, Z) - g((\nabla_{Y}S)X, Z)\\
& =  \eta(X) g(JY,Z) - \eta(Y) g(JX,Z)  - 2 \eta(Z) g(JX,Y)  \\
& \quad + g(AY,Z)g(AX, N) - g(AX,Z)g(AY, N)  \\
& \quad + \eta(AX) g(JAY,Z)-  \eta(AY) g(JAX,Z)
\end{split}
\end{equation}
where $X$, $Y$, $Z$ and $W$ are tangent vector fields of $M$.

\vskip 3pt

Note that $JX=\phi X + \eta(X)N$ and $JN=-\xi$, where $\phi X$ is the tangential component of $JX$ and $N$ is a (local) unit normal vector field of $M$. The tangent bundle $TM$ of $M$ splits orthogonally into  $TM = {\mathcal C} \oplus {\mathbb R}\xi$, where ${\mathcal C} = \mathrm{ker}\,\eta$ is the maximal complex subbundle of $TM$. The structure tensor field $\phi$ restricted to ${\mathcal C}$ coincides with the complex structure $J$ restricted to~${\mathcal C}$, and $\phi \xi = 0$. Moreover, since the complex quadric $Q^{m}$ has also a real structure~$A$, we decompose $AX$ into its tangential and normal components for a fixed $A \in \mathfrak A_{[z]}$ and $X \in T_{[z]}M$:
\begin{equation}\label{AX}
AX=BX + \rho(X)N
\end{equation}
where $BX$ is the tangential component of $AX$ and
\begin{equation*}
\rho(X)=g(AX, N)=g(X, AN)=g(X, AJ\xi) = g(JX, A \xi).
\end{equation*}
From these notations, the equations~\eqref{eq: 2.1} and \eqref{codazzi equation} can be written as
\begin{equation*}
\begin{split}
& R(X,Y)Z - g(SY,Z)SX + g(SX,Z)SY \\
& =  g(Y,Z)X - g(X,Z)Y + g(JY,Z) \phi X - g(JX,Z) \phi Y - 2g(JX,Y) \phi Z \\
& \ \  + g(AY,Z) BX - g(AX,Z)BY + g(JAY,Z) \phi BX  \\
& \ \  -  g(JAY,Z) \rho(X) \xi - g(JAX,Z) \phi BY +  g(JAX,Z)\rho(Y) \xi  \\
\end{split}
\end{equation*}
and
\begin{equation*}
\begin{split}
& (\nabla_{X}S)Y - (\nabla_{Y}S)X \\
& =  \eta(X) \phi Y - \eta(Y) \phi X  - 2 g(JX,Y)\xi  \\
& \quad + g(AX, N)BY - g(AY, N)BX  + \eta(AX) \phi BY \\
& \quad  -  \eta(AX)\rho(Y) \xi -  \eta(AY) \phi BX +  \eta(AY)\rho(X) \xi,
\end{split}
\end{equation*}
which are called the equations of Gauss and Codazzi, respectively. Moreover, from~\eqref{eq: 2.1} the Ricci tensor $\mathrm{Ric}$ of $M$ is given by
\begin{equation}\label{ricci}
\begin{split}
\mathrm{Ric}X &= (2m-1)X -3\eta(X) \xi + g(A\x, \x)BX - g(AX, N)\phi A\xi   \\
& \quad \   +g(AX, \xi)A\xi + hSX -S^{2}X,
\end{split}
\end{equation}
where $h=\mathrm{Tr}S$.

\vskip 6pt

As mentioned in section~\ref{section 2}, since the normal vector field~$N$ belongs to $T_{[z]}Q^{m}$, $[z] \in M$, we can choose $A \in {\mathfrak A}_{[z]}$ such that
\begin{equation*}
N = \cos (t) Z_1 + \sin (t) JZ_2
\end{equation*}
for some orthonormal vectors $Z_1$, $Z_2 \in V(A)$ and $0 \leq t \leq \frac{\pi}{4}$ (see Proposition~3 in~\cite{R}). Note that $t$ is a function on $M$. If $t=0$, then $N=Z_{1}\in V(A)$, therefore we see that $N$ becomes the $\mathfrak A$-principal singular tangent vector field. On the other hand, if $t=\frac{\pi}{4}$, then $N= \frac{1}{\sqrt{2}}(Z_{1}+JZ_{2})$. That is, $N$ is to be the $\mathfrak A$-isotropic singular tangent vector field. In addition, since $\xi = -JN$, we have
\begin{equation}\label{AN, Axi}
\begin{cases}
\xi  =  \sin (t) Z_2 - \cos (t) JZ_1, \\
AN  =  \cos (t) Z_1 - \sin (t) JZ_2,  \\
A\xi  =  \sin (t) Z_2 + \cos (t) JZ_1.
\end{cases}
\end{equation}
This implies $g(\xi,AN) = 0$ and $g(A\xi, \xi) = -g(AN, N)=-\cos ( 2t)$ on $M$. At each point $[z] \in M$ we define the maximal ${\mathfrak A}$-invariant subspace of $T_{[z]} M$, $[z] \in M$ as follows:
\begin{equation*}
{\mathcal Q}_{[z]} = \{X \in T_{[z]}M \mid AX \in T_{[z]}M\ {\rm for\ all}\ A \in {\mathfrak A}_{[z]}\}.
\end{equation*}
It is known if $N_{[z]}$ is ${\mathfrak A}$-principal, then ${\mathcal Q}_{[z]} = {\mathcal C}_{[z]}$ (see \cite{S1}).

\vskip 6pt

We now assume that $M$ is a Hopf hypersurface in the complex quadric~$Q^{m}$. Then the shape operator~$S$ of~$M$ in $Q^m$ satisfies $S\xi = \alpha \xi$ with the Reeb function $\alpha = g(S\xi,\xi)$ on $M$. By virtue of the Codazzi equation, we obtain the following lemma.
\begin{lemma}[\cite{SDGA}] \label{lemma Hopf}
Let $M$ be a Hopf hypersurface in $Q^m$, $m \geq 3$. Then we obtain
\begin{equation}\label{eq: 3.2}
 X \alpha = (\xi \alpha) \eta(X)  + 2g(A\xi,\xi)g(X,AN)
\end{equation}
and
\begin{equation}\label{eq: 3.1}
\begin{split}
& 2g(S \phi SX,Y) - \alpha g((\phi S + S\phi)X,Y) - 2 g(\phi X,Y) \\
& \quad  + g(X,AN)g(Y,A\xi) - g(Y,AN)g(X,A\xi)\\
& \quad  - g(X,A\xi)g(JY,A\xi) + g(Y,A\xi)g(JX,A\xi)\\
& \quad  - 2g(X,AN)g(\xi,A\xi)\eta(Y) + 2g(Y,AN)g(\xi,A\xi)\eta(X)=0
\end{split}
\end{equation}
for any tangent vector fields $X$ and $Y$ on $M$.
\end{lemma}

\begin{remark}\label{remark 3.3}
{\rm By virtue of \eqref{eq: 3.2} we know {\it if $M$ has vanishing geodesic Reeb flow (or constant Reeb curvature, respectively), then the normal vector $N$ is singular}. In fact, under this assumption \eqref{eq: 3.2} becomes $g(A \xi , \xi) g(X, AN)=0$ for any tangent vector field $X$ on $M$. Since $g(A\xi, \xi)=-\cos (2t)$, the case of $g(A\xi,\xi)=0$ implies that $N$ is $\mathfrak A$-isotropic. Besides, if $g(A\xi, \xi)\neq 0$, that is, g(AN, X)=0 for all $X \in TM$, then
$$
AN= \sum_{i=1}^{2m}g(AN, e_{i})e_{i} + g(AN, N)N = g(AN,N)N,
$$
which implies that $N=A^{2}N=g(AN, N)AN$. Taking an inner product with $N$, it follows $g(AN, N)=\pm 1$. Since $g(AN, N)=\cos (2t)$ where $t \in [0, \frac{\pi}{4})$, we obtain $AN=N$. Hence $N$ should be $\mathfrak A$-principal.}
\end{remark}

\begin{lemma}[\rm \cite{S1}]\label{lemma 3.2}
Let $M$ be a Hopf hypersurface in $Q^m$ such that the normal vector
field $N$ is ${\mathfrak A}$-principal everywhere. Then $\alpha$ is
constant. Moreover, if $X \in {\mathcal C}$ is a principal curvature
vector of $M$ with principal curvature $\lambda$, then $2\lambda \neq
\alpha$ and its corresponding vector~$\phi X$ is a principal curvature vector of $M$ with
principal curvature $\frac{\alpha\lambda + 2}{2{\lambda}-{\alpha}}$.
\end{lemma}
\begin{lemma}[\rm \cite{S1}]\label{lemma 3.3}
Let $M$ be a Hopf hypersurface in $Q^m$, $m \geq 3$, such that the
normal vector field $N$ is ${\mathfrak A}$-isotropic everywhere.
Then $\alpha$ is constant.
\end{lemma}
\noindent If the normal vector $N$ is $\mathfrak A$-isotropic, then we obtain
\begin{equation*}
g(A \x, N)=g(A\x, \x)=g(AN,N)=0
\end{equation*}
from \eqref{AN, Axi} and the notation of $N$. Taking the covariant derivative of $g(AN,N)=0$ along the direction of any $X \in T_{[z]}M$, $[z] \in M$, it becomes
\begin{equation*}
\begin{split}
0=X \big(g (AN, N) \big) & = g\big({\bar{\nabla}}_X(AN),N\big)+g \big(AN, {\bar\nabla}_XN \big)\\
&=g\big((\bar{\nabla}_{X}A) N + A (\bar{\nabla}_{X}N), N \big) +g \big(AN,\bar{\nabla}_XN\big)\\
&=g\big(q(X)JAN-ASX,N \big)-g\big(AN,SX \big)\\
&=-2g\big(ASX,N\big),
\end{split}
\end{equation*}
where we have used the covariant derivative of the complex structure~$A$, that is, $(\bar \nabla_{X}A)Y = q(X) JAY$ and the formula of Weingarten. Then the above formula gives $SAN=0$, because~$AN$ becomes a tangent vector field on $M$ for $\mathfrak A$-isotropic unit normal vector field~$N$.

\vskip 3pt

On the other hand, by differentiating $g(A{\xi},N)=0$ and using the formula of Gauss, we have:
\begin{equation*}
\begin{split}
0&=g\big({\bar{\nabla}}_X(A{\xi}),N \big)+g\big(A{\xi},{\bar{\nabla}}_XN \big)\\
&=g\big(({\bar{\nabla}}_X A) {\xi} + A ({\bar{\nabla}}_X \xi),N \big) + g\big(A{\xi},{\bar{\nabla}}_XN \big)\\
&=g\big(({\bar{\nabla}}_X A) {\xi}, N\big) + g\big( \nabla_{X}\xi + \sigma(X, \xi), AN \big) + g\big(A{\xi},{\bar{\nabla}}_XN \big)\\
&=g\big(q(X)JA{\xi}, N\big) + g \big ({\phi}SX+g(SX,{\xi})N, AN\big)-g\big(A{\xi}, SX\big)\\
&=-2g(A{\xi},SX),
\end{split}
\end{equation*}
where $\sigma$ is the second fundamental form of $M$ and ${\phi}AN=JAN=-AJN=A{\xi}$. By $g(A\xi, N)=0$, the vector field~$A\xi$ becomes a tangent vector field on $M$ with $\mathfrak A$-isotropic unit normal vector field~$N$. Then the above formula gives $SA{\xi}=0$.

\vskip 3pt

Moreover, when the normal vector $N$ is $\mathfrak A$-isotropic, the tangent vector space~$T_{[z]}M$, $[z] \in M$, is decomposed
\begin{equation*}
T_{[z]}M=[{\xi}] \oplus [A{\xi},AN] \oplus {\mathcal Q},
\end{equation*}
where ${\mathcal C} \ominus {\mathcal Q}={\mathcal Q}^{\bot}=\text{Span}[A{\xi},AN]$. From the equation~\eqref{eq: 3.1}, we obtain
\begin{equation*}
(2\lambda  - \alpha) S \phi X  = (\alpha  \lambda + 2) \phi X
\end{equation*}
for some principal curvature vector $X \in \mathcal Q \subset T_{[z]}M$ such that $S X = \lambda X$. If $2\lambda - \alpha =0$ (i.e. $\lambda = \frac{\alpha}{2}$), then $\alpha  \lambda + 2 = \frac{\alpha^{2}+4}{2} =0$, which makes a contradiction. Hence we obtain:
\begin{lemma}\label{lemma 3.5}
Let $M$ be a Hopf hypersurface in $Q^m$ such that the normal vector
field $N$ is ${\mathfrak A}$-isotropic. Then $S A \xi = 0$ and $SAN =0$. Moreover, if $X \in {\mathcal Q}$ is a principal curvature
vector of $M$ with principal curvature $\lambda$, then $2\lambda \neq
\alpha$ and its corresponding vector~$\phi X$ is a principal curvature vector of $M$ with
principal curvature $\frac{\alpha\lambda + 2}{2{\lambda}-{\alpha}}$.
\end{lemma}

On the other hand, from the property of $g(A\xi, N)=0$ on a real hypersurface~$M$ in $Q^{m}$ we see that the non-zero vector field~$A\xi$ is tangent to $M$. Hence by Gauss formula it induces
\begin{equation*}
\begin{split}
\N_{X}(A\x)& = {\bar \N}_{X}(A\x) - \sigma(X, A\x)  \\
& =q(X) JA\x + A(\N_{X}\x) + g(SX, \xi) AN - g(SX, A \x)N
\end{split}
\end{equation*}
for any $X \in TM$. From $AN=AJ\xi=-JA\xi$ and $JA\xi = \phi A\xi + \eta(A\xi) N$, it gives us
\begin{equation}\label{eq: A xi}
\left \{
\begin{array}{l}
\mathrm{Tangent\ Part:} \ \  \nabla_{X}(A \xi)  = q(X) \phi A\xi + B \phi SX - g(SX, \xi) \phi A \xi   \\
 \\
\mathrm{Normal \ Part:} \ \  q(X) g(A\xi, \xi)  = - g(AN, \N_{X}\x) + g(SX, \xi) g(A\xi, \xi) + g (SX, A\xi)
\end{array}
\right.
\end{equation}
In particular, if $M$ is Hopf, then the second equation in \eqref{eq: A xi} becomes
\begin{equation}\label{eq: 3.5}
q(\x) g(A\xi, \xi) = 2\alpha g(A\xi, \xi).
\end{equation}

\vskip 17pt

\section{Proof of Theorem~1 \\
- Reeb parallel structure Jacobi operator with $\mathfrak A$-principal normal -}\label{section 4}
\setcounter{equation}{0}
\renewcommand{\theequation}{4.\arabic{equation}}
\vspace{0.13in}

Let us $M$ be a real hypersurface in the complex quadric~$Q^{m}$, $m \geq 3$, with Reeb parallel structure Jacobi operator, that is,
\begin{equation*}\label{Reeb parallel structure Jacobi}
(\nabla_{\xi} {R}_{\xi})Y = 0
\tag{*}
\end{equation*}
for all tangent vector fields $Y$ of $M$.

\vskip 6pt

As mentioned in section~\ref{section 1} the structure Jacobi operator~${R}_{\xi}\in \mathrm{End}(TM)$ with respect to the unit tangent vector field~$\xi \in TM$ is induced from the curvature tensor~${R}$ of $M$ given in section~\ref{section 3} as follows: for any tangent vector fields $Y$, $Z \in TM$
\begin{equation*}\label{eq: 4.1}
\begin{split}
g({R}_{\xi} Y, Z) &= g({R}(Y, \xi)\xi, Z) \\
                  &= g(Y,Z) - \e(Y) \e(Z)  + g(A\x,\x)g(AY, Z) - g(Y,A\x)g(A\x, Z)\\
                  & \ \  - g(AY,N)g(AN,Z) + \alpha g(SY,Z) - \alpha^{2}\e(Y)\e(Z), \\
               \end{split}
\end{equation*}
where we have used $J\x=N$, $JA=-AJ$, and $g(A\xi, N)=0$.

\vskip 3pt

\begin{rem}
\rm For any tangent vector field $X$ of $M$ the vector field $AX$ belongs to $TQ^{m}$, that is, $AX=BX + \rho(X)N \in TM \oplus (TM)^{\bot}=TQ^{m}$. Therefore, from \eqref{eq: 4.1} the structure Jacobi operator of $M$ is given by
\begin{equation}\label{structure Jacobi operator}
\begin{split}
{R}_{\xi} Y &= Y - \e(Y) \x  + g(A\x,\x)BY - g(A\x, Y)A\x \\
            & \ \  - g(\phi A \xi, Y)\phi A\x + \alpha SY - \alpha^{2}\e(Y)\x.
\end{split}
\end{equation}
Here we have used that $A\xi =B \xi \in TM$ (i.e. $\rho (\xi) = g(AN, \x) =0$) and $AN=AJ\x = -JA\x = -\phi A\x -\eta(A\xi) N$.
\end{rem}

\noindent Taking the covariant derivative of \eqref{structure Jacobi operator} along the direction of $X \in TM$, then we have
\begin{equation*}
\begin{split}
& (\nabla_{X}{R}_{\x})Y \\
& \quad = -g(Y, \nabla_{X}\xi) \xi - \eta(Y) \nabla_{X}\xi + g\big(\nabla_{X}(A\xi), \xi \big) BY + g(A\x, \nabla_{X}\xi)BY \\
& \quad \ \ \ + g(A\x, \x) (\nabla_{X}B)Y - g\big(\nabla_{X}(A\xi), Y\big) A\xi - g(A\xi, Y) \nabla_{X}(A\xi) \\
&\quad \ \ \  - g\big((\nabla_{X}\phi)A\xi, Y\big) \phi A \xi + g\big(\nabla_{X}(A\xi), \phi Y\big) \phi A \xi \\
&\quad \ \ \  -g(\phi A \xi, Y) (\nabla_{X} \phi) A \x - g(\phi A \xi, Y) \phi \big(\nabla_{X}(A\xi)\big)\\
&\quad \ \ \   + (X \alpha) SY + \alpha (\nabla_{X}S)Y - 2\alpha (X\alpha) \e(Y) \xi  - \alpha^{2}g(Y, \nabla_{X} \xi) \xi - \alpha^{2} \e(Y) \nabla_{X}\xi\\
& \quad = -g(Y, \phi SX) \xi - \eta(Y) \phi SX + g(B \phi SX, \x) BY + g(A\x, \phi SX)BY \\
& \quad \ \ \ + g(A\x, \x)\Big\{q(X) JAY + g(SX, Y)AN - q(X) g(AY, \xi)N  \Big\} \\
& \quad \ \ \ + g(A\x, \x)\Big\{g(SX, Y) g(A\xi, \xi) N + g(AN, Y)SX \Big\} \\
&\quad \ \ \  - \Big \{ \big (q(X)-\alpha \eta(X)\big ) g(\phi A\xi, Y) + g(B \phi SX, Y)\Big \} A\xi\\
&\quad \ \ \  -g(A\x,Y)  \Big \{ \big(q(X)- \alpha \eta(X) \big) \phi A \xi + B \phi SX \Big \}\\
&\quad \ \ \  - \Big \{ g(A\xi, \xi) g(SX, Y) - g(SX, A\xi) \eta(Y) \Big \}  \phi A \xi + \big(q(X) - \alpha \eta(X)\big) g(A\xi, Y) \phi A \xi \\
&\quad \ \ \   - \Big \{ \big(q(X) - \alpha \eta(X)\big) g(A\xi, \xi) \eta(Y) - g(B \phi SX, \phi Y)  \Big \} \phi A \xi\\
&\quad \ \ \  -g(\phi A \xi, Y) \Big\{ \eta(A\xi) SX - g(SX, A \xi) \xi \Big\} \\
&\quad \ \ \ + g(\phi A \xi, Y) \Big \{ \big(q(X)-\alpha \eta(X)\big) A\xi - g(A\xi, \xi)\big(q(X)-\alpha \eta(X)\big) \xi - \phi B \phi SX \Big\} \\
&\quad \ \ \   + (X \alpha) SY + \alpha (\nabla_{X}S)Y - 2\alpha (X\alpha) \e(Y) \xi  - \alpha^{2}g(Y, \phi SX) \xi - \alpha^{2} \e(Y) \phi SX,
\end{split}
\end{equation*}
where we have used \eqref{eq: A xi} and
\begin{equation*}
\begin{split}
(\nabla_{X} B)Y & = \nabla_{X}(BY) - B (\nabla_{X}Y) \\
& = {\bar \N}_{X}(BY) - \sigma(X, BY) - B (\nabla_{X}Y)\\
& = {\bar \N}_{X}\big(AY - g(AY, N)N \big) - g(SX, BY)N - B (\nabla_{X}Y) \\
& = ({\bar \N}_{X}A)Y + A\big({\bar \N}_{X}Y \big) - g\big({\bar \N}_{X}(AY), N\big)N - g(AY, {\bar \N}_{X}N)N \\
&\quad \ \   - g(AY, N) {\bar \N}_{X}N - g(SX, BY)N - B (\nabla_{X}Y) \\
& = q(X) JAY + A(\N _{X}Y) + g(SX, Y) AN - q(X) g(JAY, N)N  \\
& \quad \ \ -g(\N _{X}Y, AN) N - g(SX, Y) g(AN, N) N  \\
& \quad \ \ + g(AY, SX)N +g(AY, N) SX - g(SX, BY)N - B (\nabla_{X}Y)\\
& = q(X) JAY  + g(SX, Y) AN - q(X) g(AY, \xi)N \\
& \quad \ \ + g(SX, Y) g(A\xi, \xi) N  +g(AY, N) SX.
\end{split}
\end{equation*}
Since $M$ is a Hopf real hypersurface in $Q^{m}$ with Reeb parallel structure Jacobi operator, it yields that
\begin{equation}\label{eq: Reeb parallel structure Jacobi}
\begin{split}
& g(A\x, \x)\Big\{q(\xi) JAY + \alpha \eta(Y)AN - q(\xi) g(AY, \xi)N  \Big\} \\
\end{split}
\end{equation}
\begin{equation*}
\begin{split}
& \quad + g(A\x, \x)\Big\{\alpha \eta(Y) g(A\xi, \xi) N + \alpha g(AN, Y)\xi \Big\} \\
&\quad   -\big(q(\xi) - \alpha \big) g(A\xi, \xi) \eta(Y) \phi A \xi
 -g(\phi A \xi, Y) g(\xi, A \xi) \big(q(\xi) - \alpha \big) \xi \\
&\quad   + (\xi \alpha) SY + \alpha (\nabla_{\xi}S)Y - 2\alpha (\xi \alpha) \e(Y) \xi =0.
\end{split}
\end{equation*}

\vskip 6pt


From now on, we assume that $M$ is a Hopf real hypersurface with non-vanishing geodesic Reeb flow and with Reeb parallel structure Jacobi operator in the complex quadric~$Q^{m}$, $m \geq 3$. In addition, we suppose that the normal vector field~$N$ of $M$ is $\mathfrak A$-principal. Then this assumption gives us
\begin{equation*}
AN=N \ \ \mathrm{and} \ \ A\xi = - \xi
\end{equation*}
from \eqref{AN, Axi}. So it follows that $AY \in TM$ for all $Y \in TM$, that is, $g(AY, N)=g(Y, AN)=0$. Moreover, taking the derivative to $AN=N$ with respect to the Levi-Civita connection~$\bar \nabla$ of $Q^{m}$ and using \eqref{eq: A xi}, we get
\begin{equation}\label{eq: ASY}
ASY=SY-2 \alpha (Y) \xi,
\end{equation}
together with $({\bar \nabla}_{Y}A)X = q(Y) JAX$ and ${\bar \nabla}_{Y}N = - SY$.

\vskip 3pt

From these properties, the equation~\eqref{eq: Reeb parallel structure Jacobi} can be rearranged as follows.
\begin{equation}\label{eq: 5.1}
\begin{split}
0&=({\N}_{\xi}R_{\xi})Y \\
 & =-q({\xi})JAY - q(\xi) \eta(Y) N +({\xi}{\alpha})SY+{\alpha}({\N}_{\xi}S)Y-2\alpha ({\xi}{\alpha}){\eta}(Y){\xi}
\end{split}
\end{equation}
In addition, from \eqref{eq: 3.5} we know $q({\xi})=2{\alpha}$. By Lemma~\ref{lemma 3.2} and our assumption, the Reeb curvature function $\alpha$ is non-zero constant on~$M$. So \eqref{eq: 5.1} reduces to the following
\begin{equation}\label{eq: 5.2}
({\N}_{\xi}S)Y = 2 {\phi}AY,
\end{equation}
together with $JAY=\phi AY + \eta(AY) N = \phi AY - \eta(Y)N$.

\vskip 3pt

On the other hand, by using the equation of Codazzi in section~\eqref{eq: 3.2}, we have
\begin{equation*}
\begin{split}
g\big(({\N}_{\xi} S){Y}-({\N}_{Y}S)\xi,Z\big)&= g({\phi}Y,Z)-g(AY,N)g(A{\xi},Z)\\
& \quad \ \ +g(\xi,A{\xi})g(JAY,Z) + g({\xi},AY)g(AN,Z)\\
& = g({\phi}Y,Z)-g({\phi}AY,Z).
\end{split}
\end{equation*}
Since $M$ is Hopf and Lemma~\ref{lemma 3.2}, it leads to
\begin{equation*}
\begin{split}
({\N}_{\xi}S)Y &= (\nabla_{Y}S)\xi + \phi Y - \phi AY \\
& = \alpha \phi SY - S \phi SY  + \phi Y - \phi AY
\end{split}
\end{equation*}
From this, together with \eqref{eq: 5.2}, it follows that
\begin{equation}\label{eq: 5.3}
\alpha \phi SY - S \phi SY  + \phi Y = 3 \phi AY.
\end{equation}
By virtue of Lemma~\ref{lemma Hopf}, for the $\mathfrak A$-principal unit normal vector field, we obtain
\begin{equation}\label{eq: 5.4}
2 S{\phi}SY= \alpha (S{\phi}+{\phi}S)Y + 2{\phi}Y.
\end{equation}
Therefore, \eqref{eq: 5.3} can be written as
\begin{equation}\label{eq: 5.5}
\alpha ({\phi}S-S{\phi})Y=6 {\phi}AY.
\end{equation}
Inserting $Y=SX$ for $X \in \mathcal C$ into \eqref{eq: 5.5} and taking the structure tensor~$\phi$ leads to
\begin{equation*}
\alpha S^{2}X  + \alpha \phi S{\phi} SX = 6 ASX,
\end{equation*}
where $\mathcal C =\mathrm{ker}\,\eta$ denotes the maximal complex subbundle of $TM$, which is defined by a distribution $\mathcal C = \{ X \in T_{[z]}M\,|\, \eta(X)=0\}$ in $T_{[z]}M$, $[z] \in M$. By using \eqref{eq: ASY} and~\eqref{eq: 5.4} this equation gives us
\begin{equation}\label{e86}
\alpha^{2} \phi S{\phi} X = -2 \alpha S^{2}X   + \alpha^{2} SX + 2 \alpha X + 12 SX
\end{equation}
for all $X \in \mathcal C$.

\vskip 6pt

On the other hand, in this subsection we have assumed that the normal vector field~$N$ of $M$ is $\mathfrak A$-principal. It follows that $AY \in TM$ for all $Y \in TM$. From this, the anti-commuting property with respect to $J$ and $A$ implies $\phi A X = -A \phi X$. Hence \eqref{eq: 5.5} can be expressed as
\begin{equation}\label{eq: 5.10}
\alpha ({\phi}S-S{\phi})Y= -6 A \phi Y.
\end{equation}
Putting $Y=\phi X$ into \eqref{eq: 5.10}, it gives
$$
\alpha {\phi}S \phi X = - \alpha SX + 6 AX
$$
for all $X \in \mathcal C$. Inserting this into \eqref{e86} gives
\begin{equation}\label{e88}
3 \alpha AX + \alpha S^{2}X - \alpha^{2} SX - \alpha X - 6 SX =0.
\end{equation}
Taking the complex conjugate~$A$ to \eqref{e88} and using \eqref{eq: ASY} again, we get
\begin{equation}\label{e89}
3 \alpha X + \alpha S^{2}X - \alpha^{2} SX - \alpha AX - 6 SX =0,
\end{equation}
for all $X \in \mathcal C$. Summing up \eqref{e88} and \eqref{e89}, gives $AX=X$ for all $X \in \mathcal C$. This gives a contradiction. In fact, it is well known that the trace of the real structure~$A$ is zero, that is, $\mathrm{Tr} A =0$ (see Lemma~1 in \cite{BS}). For an orthonormal basis $\{ \, e_{1}, e_{2} \cdots, e_{2m-2}, e_{2m-1}=\xi, e_{2m}=N \,\}$ for $T{Q^m}$, where $ e_{j} \in \mathcal C$ $(j = 1, 2, \cdots, 2m-2)$, the trace of $A$ is given by
\begin{equation*}
\begin{split}
\mathrm{Tr}A & = \sum_{i=1}^{2m} g(A e_{i}, e_{i}) \\
&= g(AN, N) + g(A \xi, \xi) + \sum_{i=1}^{2m-2} g(A e_{i}, e_{i}) \\
&= 2m-2.
\end{split}
\end{equation*}
It implies that $m=1$. But we now consider for the case $m \geq 3$.

\vskip 3pt

Consequently, this completes the proof that {\it there does not exists a Hopf real hypersurface $(\alpha \neq 0)$ in complex quadrics~${Q^{m}}$, $m \geq 3$, with Reeb parallel structure Jacobi operator and $\mathfrak A$-principal normal vector field}.

\vskip 17pt

\section{Proof of Theorem~2 \\
- Reeb parallel structure Jacobi operator with $\mathfrak A$-isotropic normal -}\label{section 5}
\setcounter{equation}{0}
\renewcommand{\theequation}{5.\arabic{equation}}
\vspace{0.13in}

In this section, we assume that the unit normal vector field $N$ is $\mathfrak A$-isotropic and~$M$ is a real hypersurface in complex quadric~$Q^{m}$ with non-vanishing geodesic Reeb flow and with Reeb parallel structure Jacobi operator. Then the normal vector field $N$ can be written as
$$N=\frac{1}{\sqrt 2}(Z_1+JZ_2)$$
for some orthonormal vectors $Z_1$, $Z_2 \in V(A)$, where $V(A)$ denotes a $(+1)$-eigenspace of the complex conjugation~$A \in~{\mathfrak A}$. Then it follows that
$$
AN=\frac{1}{\sqrt 2}(Z_1-JZ_2),\ AJN=-\frac{1}{\sqrt 2}(JZ_1+Z_2),\  \text{and}\  JN=\frac{1}{\sqrt 2}(JZ_1-Z_2).
$$
Then it gives that
$$
g({\xi},A{\xi})=g(JN,AJN)=0,\  g({\xi},AN)=0\ \text{and}\ g(AN,N)=0,
$$
which means that both vector fields $AN$ and $A\xi$ are tangent to $M$. From this and Lemma~\ref{lemma 3.3}, we see that the shape operator~$S$ of $M$ becomes to be Reeb parallel, that is, $(\nabla_{\xi}S)Y=0$ for all tangent vector field $Y$ on $M$.

\vskip 3pt

On the other hand, from the Codazzi equation~\eqref{codazzi equation} we obtain
\begin{equation*}
\begin{split}
(\nabla_{\x}S)Y & = (\nabla_{Y}S)\xi + \phi Y - g(AY, N) A\x + g(A\x, Y) AN \\
& = (Y\al)\xi + \alpha \phi SY - S\phi SY + \phi Y +g(A\x, Y) AN -g(AN, Y) A\x \\
& = \frac{\al}{2} (\phi S - S \phi)Y,
\end{split}
\end{equation*}
where the third equality holds from Lemmas~\ref{lemma Hopf} and \ref{lemma 3.3}. From this and $M$ has non-vanishing geodesic Reeb flow, we see that $M$ has isometric Reeb flow, that is, $S\phi = \phi S$.

\vskip 6pt

\noindent Consequently, we obtain:
\begin{proposition}\label{proposition 4.1}
Let $M$ be a real hypersurface with non-vanishing geodesic Reeb flow in the complex quadrics ${Q^{m}}$, $m \geq 3$. If the unit normal vector field~$N$ of $M$ is $\mathfrak A$-isotropic and the structure Jacobi operator~$R_{\xi}$ of $M$ is Reeb parallel, then the shape operator~$S$ of $M$ satisfies the property of Reeb parallelism. Moreover, it means that the Reeb flow of $M$ is isometric.
\end{proposition}

\begin{thm B}
Let $M$ be a real hypersurface of the complex quadric~$Q^{m}$, $m \geq 3$. The Reeb flow on $M$ is isometric if and only if $m$ is even, say $m=2k$, and $M$ is an open part of a tube around a totally geodesic $\mathbb C P^{k}$ in $Q^{2k}$.
\end{thm B}

\noindent Then by virtue of Theorem~$\rm B$, we assert: {\it if $M$ is a real hypersurface in $Q^{m}$, $m \geq 3$, with the assumptions given in Proposition~\ref{proposition 4.1}, then $M$ is an open part of~$(\mathcal T_{A})$}. Here the model space $(\mathcal T_{A})$ is a tube over a totally geodesic complex projective space $\mathbb CP^{k}$ in $Q^{2k}$, $m=2k$.

\vskip 6pt

From now on, let us check the converse problem, that is, the model space $(\mathcal T_{A})$ satisfies the all assumptions stated in Proposition~\ref{proposition 4.1}. In order to do this, we first introduce one proposition given in~\cite{S1}.
\begin{pro A}
Let $(\mathcal T_{A})$ be the tube of radius $0 < r < \frac{\pi}{2}$ around the totally geodesic $\mathbb C P^{k}$ in $Q^{2 k}$. Then the following statements hold:
\begin{enumerate}[\rm (i)]
\item {$(\mathcal T_{A})$ is a Hopf hypersurface.}
\item {Every unit normal vector $N$ of $(\mathcal T_{A})$ is $\mathfrak A$-isotropic and therefore can be written in the form $N=(Z_{1} + J Z_{2})/\sqrt{2}$ with some orthonormal vectors $Z_{1}$, $Z_{2} \in V(A)$ and $A \in \mathfrak A$.}
\item {$(\mathcal T_{A})$ has four distinct constant principal curvatures and the property that the shape operator leaves invariant the maximal complex subbundle $\mathcal C$ of $T(\mathcal T_{A})$ are $J$-invariant. The principal curvatures and corresponding principal curvature spaces of $(\mathcal T_{A})$ are as follows.
\begin{center}
\begin{tabular}{c|c|c}
\hline
\mbox{principal curvature} & \mbox{eigenspace}  & \mbox{multiplicity}  \\
\hline
$0$ & ${\mathbb C} (JZ_{1}+Z_{2}) $ & $2$ \\
$-\tan(r)$ & $W_{1}$ & $2k-2$\\
$\cot(r)$ & $W_{2}$ & $2k-2$ \\
$2\cot(2r)$ & ${\mathbb R JN}$ & $1$ \\
\hline
\end{tabular}
\end{center}}
\item {Each of the two focal sets of $(\mathcal T_{A})$ is a totally geodesic $\mathbb C P^{k} \subset Q^{2k}$.}
\item {$S \phi = \phi S$ (isometric Reeb flow). }
\item {$(\mathcal T_{A})$ is a homogeneous hypersurface of $Q^{2k}$. More precisely, it is an orbit of the $U_{k+1}$-action on $Q^{2k}$ isomorphic to $U_{k+1}/U_{k-1}U_{1}$, an $S^{2k-1}$-bundle over $\mathbb C P^{k}$. }
\end{enumerate}
\end{pro A}

\noindent By virtue of $\rm (i)$ and $\rm (ii)$ in Proposition~$\rm A$, $(\mathcal T_{A})$ is a Hopf real hypersurface with $\mathfrak A$-isotropic normal vector~$N$ in $Q^{m}$. Moreover, the structure Jacobi operator~$R_{\xi}$ of~$(\mathcal T_{A})$ should be Reeb parallel, because of $\alpha = 2 \cot (2r) \neq 0$, $0 < r < \frac{\pi}{2}$.

\vskip 17pt


\end{document}